\documentclass[12pt,a4paper]{article}
\usepackage{amsfonts,amsmath,amssymb,amscd}
\usepackage{fancybox}
\usepackage[latin1]{inputenc}
\usepackage[usenames]{color}
\usepackage{graphicx}
\usepackage{subcaption}
\usepackage{float}
\usepackage{fullpage}
\usepackage{epsfig,amssymb}

\newtheorem{thm}{Theorem}[section]
\newtheorem{cor}[thm]{Corollary}

\newtheorem{prop}[thm]{Proposition}

\newtheorem{rem}[thm]{Remark}

\definecolor{dollarbill}{rgb}{0.52, 0.73, 0.4}

\newcommand{\intprod}{\mathbin{\raisebox{\depth}{\scalebox{1}[-1]{$\lnot$}}}}

\begin{document}
	
\title{\textbf{Regarding the Euler-Plateau Problem with Elastic Modulus}}

\author{A. Gruber, A. P\'ampano and M. Toda}
\date{\today}

\maketitle 

\begin{abstract}
We study equilibrium configurations for the Euler-Plateau energy with elastic modulus, which couples an energy functional of Euler-Plateau type with a total curvature term often present in models for the free energy of biomembranes. It is shown that the potential minimizers of this energy are highly dependent on the choice of physical rigidity parameters, and that the area of critical surfaces can be computed entirely from their boundary data.  When the elastic modulus does not vanish, it is shown that axially symmetric critical immersions and critical immersions of disk type are necessarily planar domains bounded by area-constrained elasticae.  The cases of topological genus zero with multiple boundary components and unrestricted genus with  control on the geodesic torsion are also discussed, and sufficient conditions are given which establish the same conclusion in these cases.\\

\end{abstract}

\section{Introduction}

The development of modern mathematics owes a good deal to the theory of \emph{minimal surfaces} in $\mathbb{R}^3$, which has been of interest to mathematicians in addition to practitioners of  other scientific fields ranging from biology to architecture.  For centuries, these objects have captivated interest due to their elegance as well as their utility, and are frequently used as idealized models for elastic membranes and other physical structures.  

Since these surfaces arise naturally as equilibrium points of the area functional, they are especially amenable to study using techniques from the calculus of variations.  In 1760, J. Lagrange raised the question of how to find the surface with least area for a given fixed boundary \cite{L}, and obtained an equilibrium condition (now called an Euler-Lagrange equation) for surfaces expressed as a graph. Some years later, J. Meusnier realized that this equation could be equivalently expressed in terms of the principal curvatures of the surface.  In particular, Meusnier showed that Lagrange's equation represents precisely the condition that the sum of the principal curvatures vanish everywhere on the surface in question.  Of course, this is now understood as the elegant condition $H\equiv 0$, representing the vanishing of the mean curvature over the surface. Almost a century after this discovery, J. Plateau demonstrated that Lagrange's original problem could in fact be physically realized. By considering soap films spanning a given fixed rigid boundary \cite{Pl}, Plateau was able to generate tangible examples of minimal surfaces,  which was a huge breakthrough for the field.  In his honor, the problem of finding a minimal surface with fixed boundary is typically referred to as \emph{Plateau's Problem}. However simple to state, Plateu's Problem was notoriously difficult to solve and remained open until 1930-1931, when the general solution was found independently by J. Douglas and T. Rad\'o \cite{D,R}.

Apart from Plateau's Problem, where the boundary is regarded as immovable and prescribed, other reasonable boundary conditions have also been considered throughout the history of minimal surfaces. One example of this is the \emph{Free Boundary Problem}, where the boundary of the surface allowed to move but also constrained to lie on a fixed supporting surface. An understanding of this problem is essential in the theory of capillarity, where the surface under examination often models a fluid membrane of negligible thickness which separates two media \cite{La,Y}.  In this case, the membrane can progress and deform along the inside of the capillary, but must stay supported on its boundary.

On the other hand, less restrictions are considered in another Plateau-type exercise known as the \emph{Thread Problem}, where the boundary of the surface is allowed to vary provided that its length is unchanged \cite{Alt}. Physically, this problem can be interpreted as searching for the soap films which span an inextensible piece of ``thread", which may bend but cannot dilate or shrink. This condition of inextensibility gives rise to a constraint on the length of the boundary, which can be understood (due to a version of Lagrange's Multipliers Principle) as an energy acting on the boundary of the surface. This yields an interesting variational problem whose energy is a linear combination of the area functional of the surface and the length functional of the boundary.  At the time of writing, this problem remains unsolved in its full generality.

Following the idea of combining the area functional with a boundary energy, L. Giomi and L. Mahadevan in 2012 investigated equilibrium configurations for the surface tension of a homogeneous membrane with elastic boundary \cite{GM}. In this setting, the energy to be minimized consists of the area functional coupled to the bending energy of its boundary.  This is precisely the so-called \emph{Euler-Plateau Problem} (see e.g. \cite{BLM,CF,GM}), which is highly general and brings together two of the oldest objects of study in differential geometry and the calculus of variations: minimal surfaces and elastic curves.  From a physical perspective, the Euler-Plateau problem is understood as finding the soap films which span a pliable loop of ``fishing line'', and models the competition between the surface tension of the film and the buckling of the line induced at the boundary.

As mentioned, understanding the Euler-Plateau problem requires understanding the theory of elastic curves (i.e. elastica), which originated at the very beginning of the calculus of variations. These curves appear when trying to determine the equilibrium shape of an ideal elastic rod bent by forces and moments acting at its ends\textemdash a problem first formulated by J. Bernoulli in 1691 \cite{B}. After mixed initial results, this problem was not considered for some time.  However, its study was revitalized by D. Bernoulli (a nephew of J. Bernoulli) in a letter to L. Euler, where Bernoulli suggested to study elastic curves as critical points of the potential energy of strain under suitable constraints.  Using this newfound variational formulation of elastica, the possible qualitative types for untwisted planar rod configurations were completely classified by L. Euler in an Appendix to his book \cite{E} of 1744, although some partial results to this end were already known to J. Bernoulli.

Returning to the topic of minimal surfaces with elastic boundary, there has recently been significant interest in generalizing the boundary energy considered in the Euler-Plateau problem.  In particular, it is known that some soap films spanning a sufficiently pliable wire reduce their potential energy by twisting out of a planar configuration, and it is advantageous to have a model which reflects this phenomenon.  Mathematically, this leads to the \emph{Kirchhoff-Plateau Problem}, where the boundary is treated as a thin elastic rod subject to both bending and twisting \cite{BF,GLF,PP}.

On the other hand, given the complexity of modern physical and biophysical models, it is also natural to consider extensions of the potential energy on the interior of the surface.  This work considers one such extension known as the \emph{Euler-Plateau Problem with Elastic Modulus}, which considers the energy functional obtained from coupling the Euler-Plateau potential with a total Gaussian curvature term.  There is physical motivation for this, as it has been observed that changes in the total Gaussian curvature occur during the formation of fusion pores in lipid membranes \cite{SK}, and can potentially be used to detect stalks during membrane fusion.  The present contribution provides some mathematical analysis of this energy, and demonstrates rigidity results for equilibrium configurations which may be useful in the eventual classification of all admissible critical surfaces.

\section{Variational Problem}

To state the relevant problem more precisely, let $\Sigma$ be a compact, connected surface with boundary and consider an immersion of $\Sigma$ in the Euclidean 3-space $\mathbb{R}^3$,
$$X:\Sigma\rightarrow\mathbb{R}^3\,.$$
It will be assumed throughout this work that $X(\Sigma)$ is an oriented surface of class at least $\mathcal{C}^2$, embedded in $\mathbb{R}^3$ with sufficiently smooth boundary $\partial\Sigma$.  Therefore, it is possible to choose a unit normal vector field $\nu$ along $\Sigma$ so that the boundary $\partial\Sigma$ becomes positively oriented with respect to this choice. As is customary, no distinction is made between the abstract surface $\Sigma$ and its image $X(\Sigma) \subset \mathbb{R}^3$ when the context is clear.

At the boundary, each connected component of $\partial\Sigma$ will be represented by a sufficiently smooth arc length parameterized curve $C:I\rightarrow\mathbb{R}^3$. Each such curve  carries an arc length parameter $s\in I=\left[0,L\right]$, where $L>0$ represents the length of $C$. Using $\left(\,\right)'$ to represent the derivative with respect to arc length, it follows that $T(s):=C'(s)$ is the unit tangent vector field along $C$. The (Frenet) curvature of $C$, denoted $\kappa$, then becomes the function $\kappa(s):=\lVert T'(s)\rVert\geq 0$.

All this considered, it is possible to define the Frenet frame along $C$. This is a triple of orthonormal vector fields $\{T,N,B\}$, where $N$ and $B$ are the unit normal and unit binormal to $C$, respectively. Note that the Frenet frame is always well defined along each boundary component $C$, which follows because each $C$ is a sufficiently smooth closed curve whose curvature can vanish only at isolated points. Moreover, the closure of $C$ also implies that both the curvature $\kappa(s)$ and the (Frenet) torsion $\tau(s)$ are periodic functions of the arc length parameter.

The discussion above implies that the Frenet equations involving the curvature $\kappa$ and torsion $\tau$ of a curve $C(s)$ can be expressed as
\begin{equation}\label{feq}
\begin{pmatrix}
T\\N\\B \end{pmatrix}'=
 \begin{pmatrix}
0 &\kappa& 0\\
-\kappa&0&\tau\\
0&-\tau&0
\end{pmatrix}\begin{pmatrix}
T\\N\\B \end{pmatrix}.
\end{equation}

For curves which lie on the surface $\Sigma$, there is another natural frame which reflects the geometry of the surface seen as an ambient environment for the curve.  Called the Darboux frame, this triple of orthonormal vector fields is given by $\{T,\nu,n\}$, where $n:=T\times\nu$ denotes the (outward pointing) vector field co-normal to the boundary $\partial\Sigma$. The derivative of this frame with respect to the arc length parameter $s$ is given by
\begin{equation}\label{ds}
\begin{pmatrix}
T\\ \nu\\ n \end{pmatrix}'=
 \begin{pmatrix}
 0 &\kappa_n& \kappa_g\\
-\kappa_n&0&-\tau_g\\
-\kappa_g&\tau_g&0
\end{pmatrix}\begin{pmatrix}
T\\ \nu\\ n \end{pmatrix},
\end{equation}
where the functions $\kappa_g(s)$, $\kappa_n(s)$ and $\tau_g(s)$ are (respectively) the geodesic curvature, the normal curvature, and the geodesic torsion of the boundary relative to the immersion $X$. 

\begin{rem}
It is intuitive to consider each function $\kappa_g,\kappa_n,\tau_g$ as the ``rate of rotation'' induced by an appropriate motion of the Darboux frame.  For example, $\kappa_g = T' \cdot n$ captures the rate at which $T$ rotates into $n$ as the frame moves positively along $\partial\Sigma$.
\end{rem}

With this definition, the classical Gauss-Bonnet Theorem is expressed in its most common form as:
\begin{equation}\label{GB}
\int_\Sigma K\,d\Sigma=\oint_{\partial\Sigma}\kappa_g\,ds+2\pi\chi(\Sigma)
\end{equation}
where $K$ is the Gaussian curvature of $\Sigma$, and $\chi(\Sigma)$ denotes its Euler-Poincar\'e characteristic.

\begin{rem}
Note that some authors define the geodesic curvature $\kappa_g$ with the opposite sign, leading to a sign difference in the Gauss-Bonnet formula.
\end{rem}

The notion of Darboux frame is readily connected to the Frenet frame from before.  To see this, denote by $\theta\in\left[-\pi,\pi\right]$ the oriented angle between $N$, the normal to the boundary component $C$, and $\nu$, the normal to the surface $\Sigma$. Thus, $\theta$ is the {\em contact angle} between $\Sigma$ and the ruled, developable surface given by $(s,t)\mapsto T(s)+tN(s)$. In this context, the Darboux frame along $\partial\Sigma$ becomes the composition of the Frenet frame and a rotation in the $(N,B)$-plane.  More precisely, the relevant rotation can be expressed in complex notation as
\begin{equation}\label{cframe}
\nu+i n=e^{i\theta}\left(N+i B\right),
\end{equation}
from which \eqref{feq} and \eqref{ds} yield the relations
\begin{eqnarray}
\kappa_g&=&\kappa\sin\theta\,,\label{kg}\\
\kappa_n&=&\kappa\cos\theta\,, \label{kn}\\
\tau_g&=&\theta'-\tau\,.\label{taug}
\end{eqnarray}
Hence, knowledge of the contact angle $\theta$ together with either frame is sufficient for the construction of the other.

The energy of a homogeneous fluid membrane bounded by an elastic curve is well approximated by the \emph{Euler-Plateau energy} \cite{GM} mentioned before, which results from adding a multiple of the surface area to the bending energy of the boundary curve.  This leads to the Euler-Plateau problem, whose solutions are minimizing configurations of this system subject to appropriate conditions on boundary behavior. As this problem couples elastic phenomena with constraints on surface area, its minimizers are known to be highly complex (c.f. \cite{GM}).  This complexity is partially due to coupling phenomena that arise from competition between the surface tension of the membrane (which leads to minimizing the area) and the elasticity of the boundary (which aims to minimize the overall deformation.)

As discussed in the Introduction, a noteworthy extension of this Euler-Plateau energy can be obtained by adding another term proportional to the total Gaussian curvature of the surface.  There is physical motivation for this \cite{SK}, as this term is suggested to be one of the major contributions to the free energy of fusion stalks in phospholipids, as well as the tendency for membranes at phase boundary to form other intermediates.  For an immersion $X:\Sigma\rightarrow\mathbb{R}^3$, this means considering the \emph{Euler-Plateau energy with elastic modulus} ($E\equiv E_{\sigma,\eta,\alpha,\beta}$),
\begin{equation}\label{energy}
E[X]:=\sigma\int_\Sigma \,d\Sigma+\eta\int_\Sigma K\,d\Sigma+\oint_{\partial\Sigma}\left(\alpha\kappa^2+\beta\right)ds\,,
\end{equation}
where $\sigma>0$, $\eta\in\mathbb{R}$, $\alpha\geq 0$ and $\beta\in\mathbb{R}$ are constants motivated by physical applications. In particular, the parameter $\sigma$ is the \emph{surface tension} \cite{Tu-Ou-Yang}, $\eta$ represents the \emph{saddle-splay elastic modulus} \cite{SK,Virga}, while the boundary parameters $\alpha$ and $\beta$ are the \emph{flexural rigidity} and the \emph{edge tension} \cite{Tu-Ou-Yang}, respectively. Briefly, $\sigma$ reflects the tendency of the bounded surface to minimize area, $\eta$ determines its potential phase changes and spontaneous curvature, $\alpha$ controls the rigidity of the elastic boundary, and $\beta$ acts as a ``Lagrange multiplier'' which enforces its inextensibility.

For convenience, it will be assumed that all connected components of the boundary $\partial\Sigma$ are made of the same material, so that the flexural rigidity $\alpha$ and the edge tension $\beta$ take the same values on all boundary components.  While this is a reasonable assumption, it need not be true in all scenarios.  In fact, there is also interest in understanding the geometry of surfaces which have a mixture of elastic and inelastic boundary components.  Such surfaces are usually referred to as having a \emph{partially elastic boundary} \cite{PP}.

To study the total energy \eqref{energy}, its first variation will now be computed.  Consider $t>0$ and a one-parameter family of variations $X: \Sigma \times (-t,t) \to \mathbb{R}^3$ of the reference immersion (also called $X$), defined by $X(\epsilon) =  X+\epsilon\,\delta X+\mathcal{O}(\epsilon^2)$ for some sufficiently smooth vector field $\delta X: \Sigma \to \mathbb{R}^3$ restricting to $\delta C: \partial\Sigma \to \mathbb{R}^3$ on the boundary. In this notation, the variational derivative of a function(al) on $\Sigma$ is simply the first-order term in its Taylor expansion around $\epsilon = 0$.  In other words, the variation of the functional $\mathcal{F}[X]$ is given by 
\[\delta\mathcal{F}[X] := \frac{d}{d\epsilon} \mathcal{F}[X(\epsilon)]\bigg|_{\epsilon = 0} .\]

Note that the variations of the immersion $X$ considered here will be assumed to have tangential as well as normal components.  In contrast to the case of closed surfaces, the use of unrestricted variations is necessary when considering surfaces with boundary.  More precisely, recall that tangential variations effectively reparameterize closed surfaces and so do not yield new information about the (parameterization-invariant) functionals usually considered in geometry.  On the other hand, surfaces with boundary are generally not invariant under tangential variations, so there is a loss of information if only normal variations are considered.  Hence, all variational derivatives in the sequel will be computed with respect to arbitrarily-directed variations.

\begin{prop} 
Let $H$ denote the mean curvature of $\Sigma$, $\kappa$ denote the Frenet curvature of $\partial\Sigma$, $(\cdot)$ denote the scalar product on $\mathbb{R}^3$, $\partial_n$ denote the derivative operator in the co-normal direction, and $\nabla\cdot$ denote the (surface) divergence operator. Then, the following equations hold true for the first variations of area, total Gaussian curvature, and bending energy:
\begin{align*}
    \delta\left(\int_\Sigma\,d\Sigma\right) &= -2\int_\Sigma H \nu\cdot \delta X\,d\Sigma+\oint_{\partial\Sigma}n\cdot \delta C\,ds\,, \\
    \delta\left(\int_\Sigma K\,d\Sigma\right) &= \oint_{\partial\Sigma}\left(\left[\tau_g'\,\nu+Kn\right]\cdot \delta C+\kappa_n\,\partial_n\left[\nu\cdot \delta X\right]\right)ds\,, \\
    \delta\left(\oint_{\partial\Sigma}\left[\alpha\kappa^2+\beta\right]ds\right) &=\oint_{\partial\Sigma}\left(2\alpha T''+\left[3\alpha\kappa^2-\beta\right]T\right)'\cdot \delta C\,ds := \oint_{\partial\Sigma} J' \cdot \delta C\,ds\,.
\end{align*}
\end{prop}
\begin{rem}
The vector field $J = 2\alpha T''+\left(3\alpha\kappa^2-\beta\right)T$ comes from the conserved Noether current associated to the invariance of elastic curves under translations in space \cite{LS}. Physically, it can be identified as the per-area contact force, which is necessarily conserved along elasticae  \cite{Jelas}. 
\end{rem} 
\textit{Proof.}  As the argument relies on standard calculations, we merely sketch the details.  For more information see e.g. \cite[Appendix A and B]{Anthony}. First, note  that the area element on $\Sigma$ varies as \[\delta\left(d\Sigma\right)=\left(-2H\nu\cdot\delta X+\nabla\cdot\left(\delta X\right)^T\right)d\Sigma\,,\]
where $(\delta X)^T$ denotes the tangential projection of $\delta X$.  The desired expression for the variation of area then follows from Stokes' Theorem and the fact that $n\cdot (\delta X)^T = n\cdot \delta C$ on $\partial\Sigma$.  On the other hand, the variation of total Gaussian curvature follows from the pointwise equation
\[ \delta K=\nabla\cdot\left( \nabla\left[\nu\cdot\delta X\right]\intprod\left[d\nu+2H\, {\rm Id}\right]\right)+2HK\nu\cdot\delta X+\nabla K\cdot \delta X\,. \]  Note that the expression $Y \intprod\, \omega := \omega(Y)$ denotes the left contraction of the vector field $Y$ with the differential one-form $\omega$. Moreover, integration by parts and the fact that $\nabla K \cdot \delta X = \nabla K \cdot (\delta X)^T$ yield
\begin{eqnarray*}
\delta\left(\int_\Sigma K\,d\Sigma\right)&=&\oint_{\partial\Sigma}\left( \nabla\left[\nu\cdot\delta X\right]\intprod\left[d\nu+2H\, {\rm Id}\right] +K\delta X\right)\cdot n\,ds
\end{eqnarray*}
from which the noted expression follows by component-wise rearrangement.  (A different approach to this calculation can be found in Appendix A.)  Finally, the variation of the bending energy follows from standard arguments involving the variation of the Frenet frame along $\partial\Sigma$ (see e.g. \cite{LS}).\hfill$\square$ 
\\

Combining the information above, the first variation formula for the Euler-Plateau energy with elastic modulus \eqref{energy} is given by
\begin{equation*}
\delta E[X]=-2\sigma\int_\Sigma H\nu\cdot \delta X\,d\Sigma+\eta\oint_{\partial\Sigma}\kappa_n\,\partial_n\left(\nu\cdot\delta X\right)ds+\oint_{\partial\Sigma}\left(J'+\left[\sigma+\eta K\right]n+\eta\tau_g'\,\nu\right)\cdot \delta C\,ds,
\end{equation*}
which directly implies the following.

\begin{thm}
The Euler-Lagrange equations for equilibria of the total energy $E[X]$ are 
\begin{eqnarray}
H&\equiv&0\,,\quad\quad\quad\text{on $\Sigma$}\,,\label{EL1}\\
\eta\kappa_n&=&0\,,\quad\quad\quad\text{on $\partial\Sigma$}\,,\label{EL2}\\
J'\cdot\nu+\eta\tau_g'&=&0\,,\quad\quad\quad\text{on $\partial\Sigma$}\,,\label{EL3}\\
J'\cdot n-\eta\tau_g^2+\sigma&=&0\,,\quad\quad\quad\text{on $\partial\Sigma$}\,.\label{EL4}
\end{eqnarray}
\end{thm} 
\begin{rem}
These Euler-Lagrange equations can also be recovered using (77)-(80) of \cite{Tu-Ou-Yang}, but they are written here in a way which is more reflective of the invariances of the problem.  As can be seen from the introduction of the Noether current $J$, this leads to a significant amount of non-obvious simplification in the boundary conditions.
\end{rem} 
\textit{Proof.}  Consider an equilibrium immersion $X$, so that $\delta E[X] = 0$.  If $\delta X$ is compactly supported, it is clear that $H\equiv 0$ must hold everywhere on $\Sigma$ and therefore the immersion $X$ must be minimal. Moreover, examining normal variations $\delta X=\psi\nu$ for some (sufficiently smooth) function $\psi: \Sigma \to \mathbb{R}$  leads to the boundary integrals
\[ 0=\oint_{\partial\Sigma}\left(J'\cdot \nu+\eta\tau_g'\right)\psi\,ds+\eta\oint_{\partial\Sigma}\kappa_n\,\partial_n\psi\,ds\,, \]
from which it follows that both $\eta\kappa_n = 0$ and $J'\cdot\nu+\eta\tau_g'=0$ hold everywhere on $\partial\Sigma$, since $\psi$ and $\partial_n\psi$ are arbitrary functions which can be prescribed on the boundary. The condition $J'\cdot n+\eta K+\sigma=0$ is then deduced in a similar manner by taking variations tangential to the immersion. Finally, recall that the Gaussian curvature along $\partial\Sigma$ is defined to be $K:=-\det(d\nu)=\kappa_n(2H-\kappa_n)-\tau_g^2$, so that the condition $\eta\kappa_n=0$ implies that $\eta K=-\eta\tau_g^2$ on $\partial\Sigma$. \hfill$\square$ 
\\

Certainly equation \eqref{EL1} implies that the equilibria of \eqref{energy} are minimal.  This is expected, as the surface tension will always encourage the interior of a critical surface to minimize its area. On the other hand, the boundary conditions \eqref{EL2}-\eqref{EL4} can be interpreted as the force and momentum equilibria for the problem (c.f. \cite{Tu-Ou-Yang}).  Note that equations \eqref{EL3} and \eqref{EL4} represent a generalization of the classical  Euler-Lagrange equations for elastic curves \cite{LS,PP}, which can also be related to the equilibrium equations for an elastic rod in the presence of an external directed force (see e.g. \cite{W}).  Moreover, it is evident from \eqref{EL2} that the boundary $\partial\Sigma$ is composed of closed asymptotic lines when $\eta\neq 0$. Note that the presence of the constant $\sigma$ in \eqref{EL4} reflects the coupling mentioned previously between the surface tension and the elasticity of the boundary.

\section{Equilibrium Configurations}

Before discussing the geometry of critical immersions for the Euler-Plateau energy with elastic modulus \eqref{energy}, it is useful to consider how this energy is affected under rescalings of a given surface immersion.  In particular, consider  $X:\Sigma\rightarrow\mathbb{R}^3$ and a rescaling $X\mapsto \lambda X$ where $\lambda > 0$. A  computation shows that this transformation scales the surface area quadratically and leaves the total Gaussian curvature invariant.  Moreover, it is straightforward to show that this rescaling also induces a linear change in the length of the boundary and an inverse linear (order $\lambda^{-1}$) change in its bending.  Consequently, it is possible to formulate an interesting relationship between the area of a critical immersion and its boundary data.

\begin{prop}\label{rescalings} Let $X:\Sigma\rightarrow\mathbb{R}^3$ be a critical immersion for the total energy $E[X]$ defined in  \eqref{energy}. Then, the following relation holds:
$$2\sigma\mathcal{A}[X]=\oint_{\partial\Sigma}\left(\alpha\kappa^2-\beta\right)ds\,,$$
where $\mathcal{A}[X]$ denotes the area of the immersion.
\end{prop}
\textit{Proof.} As argued above, after rescaling $X\mapsto\lambda X$ the total energy is given by
$$E[\lambda X]=\sigma\lambda^2\mathcal{A}[X]+\eta\int_\Sigma K\,d\Sigma+\frac{\alpha}{\lambda}\oint_{\partial\Sigma}\kappa^2\,ds+\beta\lambda\oint_{\partial\Sigma}\,ds\,.$$
Thus, differentiating the above with respect to $\lambda$ shows that at critical points ($\lambda=1$),
$$0=2\sigma\mathcal{A}[X]-\alpha\oint_{\partial\Sigma}\kappa^2\,ds+\beta\oint_{\partial\Sigma}\,ds\,,$$
which implies the conclusion. \hfill$\square$

\begin{rem}
Since the constant $\sigma$ can be obtained from equation \eqref{EL4} using only the data of the boundary, Proposition \ref{rescalings} implies that the area of a critical immersion is completely determined by the elastic energy of its boundary. See also Proposition 2.1 of \cite{PP2} for $H=0$ and $c_o \neq 0$. 
\end{rem}

To further discuss the properties of equilibrium configurations for \eqref{energy}, consider the case where $\eta\neq 0$. From \eqref{kn} and the Euler-Lagrange equation \eqref{EL2}, it is clear that the contact angle $\theta$ is constant and satisfies $\theta\equiv \pm\pi/2$. This is suggestive of the study of capillarity, where this case is distinguished even among the already restrictive case of membranes with constant contact angle \cite{P}. In particular, a bead of liquid whose normal makes a contact angle of $\pm\pi/2$ with a solid surface implies that the surface is ``perfectly wetted'' by the liquid. This means that the molecules of the liquid have perfect tendency to interact with the molecules of the solid, and are not influenced by intra-molecular interactions within the liquid itself \cite{P1,Z}. In the present case, the exceptional contact angle is a consequence of the vanishing of the normal curvature along $\partial\Sigma$, which makes the boundary a closed asymptotic line in the shared surface.

This idea has significant consequences on the energy of $E$-critical surfaces.  Using $\eta \neq 0$ and equation \eqref{2} along with the definition $J := 2\alpha T'' + \left( 3\alpha\kappa^2 - \beta \right)T$, the boundary conditions \eqref{EL3} and \eqref{EL4} can be rewritten (resp.) as
\begin{eqnarray}
4\alpha\kappa_g'\tau_g+2\alpha\kappa_g\tau_g'+\eta\tau_g'&=&0\,,\quad\quad\quad\text{on $\partial\Sigma$}\,,\label{1}\\
2\alpha\kappa_g''+\left(\alpha\kappa_g^2-2\alpha\tau_g^2-\beta\right)\kappa_g-\eta\tau_g^2+\sigma&=&0\,,\quad\quad\quad\text{on $\partial\Sigma$}\,.\label{2}
\end{eqnarray}
Note that equation (13) can be characterized as the binormal component of the Euler-Lagrange operator associated to the curvature energy representing \emph{elastic curves circular at rest},
$$\mathbf{\Theta}[C]:=\int_C \left(\left[\kappa+\mu\right]^2+\lambda\right)ds\,,$$
where $\mu:=\pm\eta/(2\alpha)$ and $\lambda:=\beta/\alpha-\mu^2$. At the same time, equation \eqref{2} is an extension of the normal component of the Euler-Lagrange operator of $\mathbf{\Theta}$. Since this extension involves the surface tension $\sigma>0$, it again illustrates the significant interaction between the surface and the boundary which takes place during the minimization of $E[X]$. 

\begin{rem}
The energy $\mathbf{\Theta}$ has also been used to study the shape of stiff rods which are circular in their undeformed state \cite{CCG}.
\end{rem}

As one of the two Euler-Lagrange equations describing critical curves for $\mathbf{\Theta}$, it follows that \eqref{1} can be integrated using a technique involving Killing vector fields along $C$ (c.f. \cite{LS}).  This yields the helpful  \emph{geodesic curvature-torsion integrable system} along each connected component $C \subset \partial\Sigma$,
\begin{equation}
\tau_g\left(2\alpha\kappa_g+\eta\right)^2=c\,.\label{fi1}
\end{equation}
In particular, the case where $c=0$ corresponds precisely to the case where the translational and rotational Noether currents associated to $\mathbf{\Theta}$ are orthogonal along $\partial\Sigma$. Further, it follows that if $c=0$ in \eqref{fi1}, then either $\tau_g\equiv 0$ or $2\alpha\kappa_g+\eta\equiv 0$ \emph{identically} on $C$. To see this, assume that $c=0$ and there exists a point $p \in C \subset \partial\Sigma$ such that $2\alpha\kappa_g+\eta\neq 0$. By continuity, there must exist a small boundary neighborhood $U$ containing $p$ on which $2\alpha\kappa_g+\eta\neq 0$, so that $\tau_g = 0$ also holds on $U$. Moreover, $\tau_g$ is real analytic as a solution to the ODE system \eqref{1}-\eqref{2} with real analytic coefficients, so this implies that $\tau_g\equiv 0$ must hold on the entirety of $C$. Of course, similar reasoning implies the conclusion $2\alpha\kappa_g + \eta \equiv 0$ on $C \subset \partial\Sigma$ when $\tau_g \neq 0$ somewhere on $C$.

\begin{rem}\label{taugnot0} If  $2\alpha\kappa_g+\eta\equiv 0$ holds identically along $\partial\Sigma$, then $\Sigma$ is a minimal surface bounded by asymptotic closed curves with constant geodesic curvature. Moreover, from \eqref{2} this is only possible when the energy parameters satisfy $8\sigma\alpha^2=\eta\left(\eta^2-4\alpha\beta\right)$.
\end{rem}

To proceed with the study of the critical points of \eqref{energy}, we first consider the case where the flexural rigidity $\alpha = 0$ at the boundary. Due to the presence of the elastic modulus term in $\eqref{energy}$, this gives an extension of the Thread Problem \cite{Alt}, which consists of searching for minimal surfaces whose boundary is an inextensible piece of ``thread''. Here, the inextensibility of $\partial\Sigma$ is enforced via the potentially nonvanishing edge tension $\beta$.  Moreover, the condition $\alpha = 0$ implies that total energy $E[X]$ is agnostic regarding any bending or twisting which occurs along the boundary. Consequently, there is the following result.

\begin{thm}\label{alpha0} Let $X:\Sigma\rightarrow\mathbb{R}^3$ be a critical immersion for $E[X]$ defined in \eqref{energy} with $\alpha=0$ and $\eta\neq 0$. Then, $\beta<0$ and the surface is a planar disk (i.e. a topological disk contained in a plane) bounded by a circle of radius $-\beta/\sigma$.
\end{thm}
\textit{Proof.} Consider an immersion $X$ critical for \eqref{energy}, with $\alpha=0$ and $\eta\neq 0$. Then, the Euler-Lagrange equations \eqref{EL1}-\eqref{EL4} are satisfied. In particular, the first integral \eqref{fi1} shows that $\eta^2\tau_g = c$ for some $c\in\mathbb{R}$ and the geodesic torsion $\tau_g$ is constant along each boundary component (possibly for different values of $c$).

Combining this with \eqref{kn}, \eqref{taug}, \eqref{EL2} and \eqref{EL3} shows that each boundary component is a Frenet helix, meaning $\kappa_g=\pm\kappa$ and $\tau_g=-\tau$ are constant.  Of course, in order for each component to ``close up'' as required, it must follow that $\tau_g=-\tau=0$ and each boundary curve is a circle. 

Now, recall that every constant mean curvature (CMC) surface admits a real analytic parameterization \cite{BL}, so the Cauchy-Kovalevskaya Theorem implies that any minimal surface \eqref{EL1} with a circular boundary component on which $\tau_g\equiv 0$ holds must be axially symmetric, i.e. a section of a plane or a catenoid (for details see Proposition 5.1 of \cite{PP2}). Since the minimal surface in the present case also satisfies $\kappa_n\equiv 0$ on $\partial\Sigma$, this implies that the surface is a planar disk. Observe that any other topologies are discarded here, since the energy parameters must be the same in all boundary components.

Finally, this information in combination with the Euler-Lagrange equation \eqref{EL4} shows that $\Sigma$ is a planar disk bounded by a circle of radius $-\beta/\sigma$. Clearly, this also implies that $\beta<0$ holds. \hfill$\square$
\\

In view of Theorem~\ref{alpha0}, it will be assumed in the sequel that the functional $E[X]$ satisfies $\alpha>0$. In addition, $2\alpha\kappa_g+\eta\neq 0$ will be assumed to hold on at least one point of the boundary component(s) under consideration (c.f. Remark~\ref{taugnot0}).

It is also useful to develop an understanding of the elasticae which bound $E$-critical domains. To that end, consider when $\tau_g=0$ holds on at least one point of a boundary component $C$. When this occurs, it follows from the first integral \eqref{fi1} that the constant of integration $c= 0$ on $C$. Since (by assumption)  $2\alpha\kappa_g+\eta\neq 0$ on at least one point of $C$, this means $\tau_g\equiv 0$ must hold on $C$ as explained above. Combining equation \eqref{taug} with the fact that $\theta\equiv\pm\pi/2$ then implies that the Frenet torsion vanishes along the boundary component $C$, i.e. $\tau\equiv 0$ holds. Therefore, $C$ is planar and contained in a suitable plane. In this case, the Euler-Lagrange equation \eqref{EL4} reduces to
\begin{equation}\label{area-cons}
2\alpha\kappa_g''+\left(\alpha\kappa_g^2-\beta\right)\kappa_g+\sigma=0\,,
\end{equation}
which is a special case of the classical second order Riccati equation with constant coefficients used to represent the Euler-Lagrange equation of \emph{area-constrained (planar) elastic curves} (see e.g. \cite{ACCG}). 

Area-constrained planar elasticae first appeared in 1884 in the work of Levy \cite{Levy}, where they were used to model thin elastic rods under a constant perpendicular force directed along their length. Note that this was in contrast with previous work on this subject, where only forces acting at the ends had been considered \cite{E}.  For this reason, such curves are sometimes referred as to \emph{elasticae under pressure}. In particular, when they are closed (as for the present case), they are known as \emph{buckled rings}. For a survey of these objects and their relation to different variational problems see e.g. \cite{W}.

Looking first for constant curvature solutions to \eqref{area-cons}, i.e. solutions where  $\kappa_g\equiv \kappa_o\in\mathbb{R}$, shows that the area-constrained elastica $C$ is necessarily a circle (note that since $C$ is closed, $\kappa_o\neq 0$). This leads to the following existence result (for some choices of the parameters we also obtain uniqueness).

\begin{prop} For fixed constants $\sigma>0$, $\alpha>0$ and $\beta\in\mathbb{R}$, there always exist area-constrained elastic circles satisfying \eqref{area-cons}. Moreover, if $27\alpha\sigma^2> 4\beta^3$ then there is only one.
\end{prop}
\textit{Proof.} The curvature of an area-constrained elastic circle is a nonzero constant $\kappa_g\equiv \kappa_o\neq 0$ which is a root of the polynomial \eqref{area-cons}
$$Q(\kappa_o)=\alpha\kappa_o^3-\beta\kappa_o+\sigma=0\,.$$
Since the limit of $Q(\kappa_o)$ when $\kappa_o\rightarrow\pm\infty$ is $\pm\infty$, respectively, there is always at least one area-constrained elastic circle. 

Moreover, after differentiating $Q(\kappa_o)$ with respect to $\kappa_o$, it is easy to see that the polynomial is non-decreasing when $\beta\leq 0$. Together with $Q(0)=\sigma>0$, this implies that there is only one negative root, and consequently there is only one area-constrained elastic circle. On the other hand, if $\beta>0$ holds, the critical points of this polynomial are a local maximum (for a negative value of $\kappa_o$) and a local minimum (for a positive value of $\kappa_o$). It is straightforward to check that the value of $Q(\kappa_o)$ at the local minimum is positive if and only if $27\alpha\sigma^2> 4\beta^3$. Hence, when this occurs, there is only one negative root of $Q(\kappa_o)$. Conversely, if $27\alpha\sigma^2=4\beta^3$ there are two roots  (one positive and another one negative), and if $27\alpha\sigma^2<4\beta^3$ there are three roots (one negative and two positive).\hfill$\square$ 
\\

Moving further, if area-constrained elastic curves $C$ with non-constant curvature are considered, then equation \eqref{area-cons} can be integrated once. Indeed, multiplying by $\kappa_g'$ yields an exact differential equation whose first integral is
\begin{equation}\label{int}
\left(\kappa_g'\right)^2=d-\frac{1}{4}\kappa_g^4+\frac{\beta}{2\alpha}\kappa_g^2-\frac{\sigma}{\alpha}\kappa_g\,,
\end{equation}
where $d\in\mathbb{R}$ is a constant of integration. From equation \eqref{int}, the expression of the curvature $\kappa_g$ can be explicitly obtained in terms of elliptic integrals (see e.g. \cite{ACCG,W} and their references). In the same references some figures are all shown.

With some understanding of critical curves, we now study the critical domains where $\tau_g=0$ holds on at least one point of the boundary component $C$. As argued above, when $2\alpha\kappa_g + \eta \neq 0$ on at least one point of $C$, equation \eqref{fi1} implies that $\tau_g\equiv 0$ holds along the entirety of $C$. Next, $\eta\neq 0$ and \eqref{EL1}, \eqref{EL2} combine to show that $\kappa_n = H = 0$ along $C$, which implies that the Gaussian curvature along $C$ is
$$K:=-{\rm det}(d\nu)=\kappa_n\left(2H-\kappa_n\right)-\tau_g^2=-\tau_g^2=0\,.$$
Therefore, the minimal surface $\Sigma$ is flat along $C$. 

\begin{rem}\label{plfl} In general, the notions of ``planar'' and ``flat'' are not equivalent, as \emph{planar} means that the object under consideration is contained in a suitable plane,  while \emph{flat} indicates that the Gaussian curvature vanishes (i.e. $K\equiv 0$) along that object. Clearly, the above computation implies that a flat boundary component in a minimal surface is planar, but the converse is not necessarily true.  On the other hand, these notions are equivalent for the boundary components of $E$-critical immersions, since the contact angle satisfies $\theta\equiv\pm\pi/2$.  
\end{rem}

To continue with the characterization of equilibrium configurations for \eqref{energy}, we need the following result concerning generic minimal immersions which are flat along a boundary component.  Since this result has several illustrative proofs, some alternatives to the one given here are sketched in Appendix B.

\begin{prop}\label{flat} Let $X:\Sigma\rightarrow\mathbb{R}^3$ be a minimal immersion of a connected surface $\Sigma$ with boundary $\partial\Sigma$. If the boundary is flat along some connected component $C\subset\partial\Sigma$, then $\Sigma$ is a planar domain.
\end{prop}
\textit{Proof.} Since the immersion $X:\Sigma\rightarrow\mathbb{R}^3$ is minimal, its image has a \emph{Weierstrass representation} (see e.g. \cite{N2}).  To elaborate, there is an analytic function $f$ and a meromorphic function $g$ so that $fg^2$ is analytic on $\Sigma$ and the image can be parameterized as 
\begin{equation}\label{Weierstrass}
X(z)=\frac{1}{2}\,\Re\left(\int_{z_o}^z\left(f\left[1-g^2\right],i f\left[1+g^2\right],2fg\right)d\omega\right).
\end{equation}

With respect to this parameterization, the Gaussian curvature $K$ of the minimal immersion $X:\Sigma\rightarrow\mathbb{R}^3$ is given by
$$K=-\frac{4\lvert dg\rvert^2}{\lvert f\rvert^2\left(1+\lvert g\rvert^2\right)^4}\,,$$
where $dg$ is a meromorphic differential one-form. It follows that the Gaussian curvature of $X$ is nonpositive, and that $K=0$ at some point if and only if $dg=0$ at that point.

Now, a meromorphic differential form is either identically zero, or its zeros are isolated. Applying this to $dg$ and using the hypothesis that $K$ (hence $dg$) vanishes along an entire boundary component, it follows that $dg\equiv 0$ on $\Sigma$.  In other words, the surface $\Sigma$ must be everywhere flat. Combining this with the minimal condition yields the conclusion that $\Sigma$ is planar. \hfill$\square$ 
\\

Using Proposition \ref{flat}, there is the following characterization of critical domains for \eqref{energy} with $\eta\neq 0$ and such that $\tau_g=0$ holds at one or more boundary points.

\begin{thm}\label{thmplanar} Let $X:\Sigma\rightarrow\mathbb{R}^3$ be a critical immersion for $E[X]$ defined in \eqref{energy} with $\eta \neq 0$ and such that $2\alpha\kappa_g + \eta \neq 0$ on at least one point of any boundary component $C\subset \partial\Sigma$. If $\tau_g=0$ holds somewhere along $C$, then the surface is a planar domain bounded by area-constrained elasticae.
\end{thm}
\textit{Proof.} Assume that $\tau_g=0$ holds on at least one point of some boundary component $C$. By previous argument, it follows that $\tau_g\equiv 0$ on $C$ and that the critical immersion of $\Sigma$ is minimal and flat along $C$. Proposition \ref{flat} then applies to conclude that the critical surface is planar.  Finally, the Euler-Lagrange equation \eqref{EL4} reduces to \eqref{area-cons} since $\tau_g\equiv 0$ along the boundary, so that $\partial\Sigma$ is composed of area-constrained planar elastic curves. \hfill$\square$
\\

Additionally, note that if a critical surface is axially symmetric, then $\tau_g\equiv 0$ holds (everywhere) along the boundary $\partial\Sigma$. Using this together with $\kappa_n\equiv 0$ along $\partial\Sigma$ (from \eqref{EL2}), the conditions of Theorem~\ref{thmplanar} are satisfied and lead directly to the following corollary.

\begin{cor} Let $X:\Sigma\rightarrow\mathbb{R}^3$ be an axially symmetric critical immersion for $E[X]$, \eqref{energy}, with $\eta\neq 0$. Then, the surface is a planar disk bounded by an area-constrained elastic circle.
\end{cor}

\section{Equilibrium Configurations of Genus Zero}

It is natural to consider the properties of $E$-critical immersions which are particular to surfaces with topological genus zero.  To that end, let $\eta\neq 0$ and $X:\Sigma\rightarrow\mathbb{R}^3$ be an immersion of a genus zero surface $\Sigma$ with boundary $\partial\Sigma$ which is critical for the total energy \eqref{energy}.  First, suppose $\Sigma\cong D$ is congruent to a topological disk.  Adapting an argument due to Nitsche \cite{N}, it is possible to establish the following result.

\begin{thm}\label{disk} Let $\eta\neq 0$ and $X:\Sigma\cong D\rightarrow\mathbb{R}^3$ be an immersion of disk type critical for the total energy $E[X]$  \eqref{energy}. Then, the surface is a planar domain bounded by an area-constrained elastic curve.
\end{thm}
\textit{Proof.} We may assume that the surface is given by a conformal immersion of the unit disk $\mathcal{D}$ in the complex plane $\mathbb{C}$. Let $z$ denote the usual complex coordinate in the disk and let $\omega:=\log z$. Although $\omega$ is not well defined, its differential $d\omega=dz/z$ is well defined in $\mathcal{D}-\{0\}$. Therefore, the fundamental forms of the immersion can be expressed in a neighborhood of $\partial\mathcal{D}$ as
\begin{equation*}
ds_X^2:=e^{\zeta}\lvert d\omega\rvert^2\,,\quad\quad\quad \mathbb{I}:=\frac{1}{2}\Re\left(\Phi \,d\omega^2\right),
\end{equation*}
where $\Phi\,d\omega^2=-\left(\mathbb{I}_{22}+i\mathbb{I}_{12}\right)d\omega^2$ is the \emph{Hopf differential} \cite{H}. Here, $\mathbb{I}_{ij}$, $i,j=1,2$ are the coefficients of the second fundamental form. Notice that $\mathbb{I}_{11}=-\mathbb{I}_{22}$ holds since the surface is minimal. Moreover, using subscript letters to denote differentiation with respect to the subscripted variable, there are the following Gauss and Codazzi equations valid for minimal immersions, 
$$\lvert\Phi\rvert^2 e^{-2\zeta}=-K\,,\quad\quad\quad \Phi_{\bar{\omega}}=0\,.$$
In particular, the second of these implies that $\Phi$ defines a complex analytic function. Additionally, it follows from $\kappa_n \equiv 0$, \eqref{EL2}, and the definition of $\Phi$ that $\Re(\Phi)\equiv 0$ on $\partial\mathcal{D}$, so that
\begin{equation}\label{Phiboundary}
\Phi\lvert_{\partial\mathcal{D}}=-i \tau_g e^\zeta.
\end{equation}

On the other hand, the transformation law for quadratic differentials yields the following relation between the Hopf differential in the $\omega$ and $z$ coordinates,
$$\Phi\, d\omega^2=\Phi \,\omega_z^2\, dz^2=\left(\frac{\Phi}{z^2}\right) dz^2=:\widetilde{\Phi}\, dz^2\,.$$
In contrast to $\Phi$, the function $\widetilde{\Phi}$ is globally defined and analytic on $\mathcal{D}$, as is $z^2\widetilde{\Phi}$. The calculation above shows that $\Phi=z^2\widetilde{\Phi}$ on $\partial\mathcal{D}$, so it follows that  $\Re(z^2\widetilde{\Phi})\equiv 0$ holds on $\partial\mathcal{D}$. However, analyticity implies that $\widetilde{\Phi}=ic/z^2$ holds on $\mathcal{D}$ for a real constant $c\in\mathbb{R}$, which is impossible unless $c=0$ and hence $\widetilde{\Phi}$ vanishes identically. It follows that if $\widetilde{\Phi}\equiv 0$ holds in $\mathcal{D}$ then every point is umbilic, which means the surface is planar.  Finally, these planar domains must satisfy the Euler-Lagrange equation \eqref{EL4} on the boundary, which is equivalent to \eqref{area-cons} since $\tau_g\equiv 0$ along $\partial\mathcal{D}$. This finishes the proof. \hfill $\square$
\\

Although the global existence of non-planar disk type critical domains is completely restricted, it is always possible to construct such domains locally by solving \emph{Bj\"orling's problem} \cite{Bjorling}, an outline of which will now be given. First, rewrite the boundary conditions \eqref{EL3} and \eqref{EL4} in terms of the Frenet curvature $\kappa(s)$ and torsion $\tau(s)$. Since $\theta\equiv \pm\pi/2$, it follows that $\kappa_g=\pm\kappa$, $\kappa_n\equiv 0$ and $\tau_g=-\tau$ hold on $\partial\Sigma$. With this, \eqref{EL3}-\eqref{EL4} become the respective equations,
\begin{eqnarray}
4\alpha\kappa'\tau+2\alpha\kappa\tau'\mp\eta\tau'&=&0\,,\quad\quad\quad\text{on $\partial\Sigma$}\,,\label{el1}\\
2\alpha\kappa''+\left(\alpha\kappa^2-2\alpha\tau^2-\beta\right)\kappa\mp \eta\tau^2\pm \sigma&=&0\,,\quad\quad\quad\text{on $\partial\Sigma$}\,.\label{el2}
\end{eqnarray}

By the Fundamental Theorem of Curves, given functions $\kappa(s)$ and $\tau(s)$ there exists a unique arc length parameterized curve, up to rigid motions, whose curvature and torsion are $\kappa(s)$ and $\tau(s)$, respectively. Let $C(s)$ be such a curve whose curvature $\kappa(s)$ and torsion $\tau(s)$ are solutions of \eqref{el1}-\eqref{el2}. Since the coefficients of these equations are real analytic, both $\kappa(s)$ and $\tau(s)$ are real analytic functions of the arc length parameter $s$. Moreover, since $C(s)$ can be found by solving \eqref{feq}, it follows that $C(s)$ is also real analytic.

Define a unit vector field $\nu(s)$ along $C(s)$ which is orthogonal to $T(s)$ and makes an angle $\theta\equiv \pm\pi/2$ with the Frenet normal $N(s)$ (i.e. from \eqref{cframe}, $\nu(s)\equiv\pm B(s)$). By analyticity, both the curve $C(s)$ and the vector field $\nu(s)$ have analytic extensions $C(z)$ and $\nu(z)$ to a simply connected domain $U \subset \mathbb{C}$  with coordinate $z=s+it$.

Next, for fixed $s = s_o$, Bj\"orling's formula (introduced by Schwarz \cite{Sch}),
$$X(z):=\Re\left(C(z)+i\int_{s_o}^z\left[C'(\omega)\times\nu(\omega)\right]d\omega\right),$$
gives a minimal surface containing the curve $C$ which has unit normal $\nu(s)$ along $C$. Finally, consider a small part of this minimal surface lying on one side of the curve $C$. Clearly, in this local domain the Euler-Lagrange equation \eqref{EL1} holds, and it follows from the choice of $\nu(s)$ that $\kappa_n\equiv 0$ along $C$. Hence \eqref{EL2} is satisfied along $C$, and the construction of the curve $C(s)$ implies that also \eqref{EL3}-\eqref{EL4} hold, so that this surface is indeed a ``local" critical domain.

\begin{figure}[h!]
\centering
\begin{subfigure}[b]{0.3\linewidth}
\includegraphics[width=\linewidth]{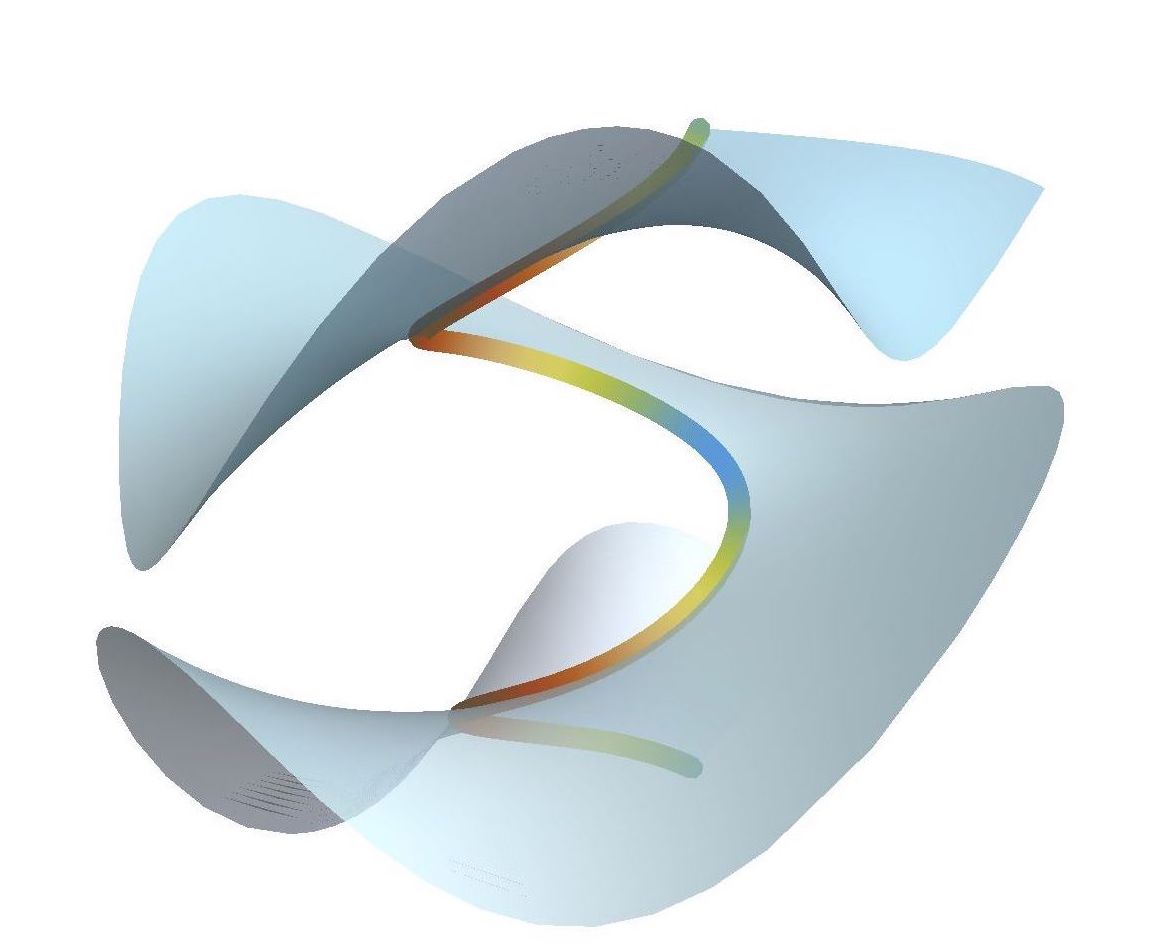}
\end{subfigure}
\quad\,\,
\begin{subfigure}[b]{0.3\linewidth}
\includegraphics[width=\linewidth]{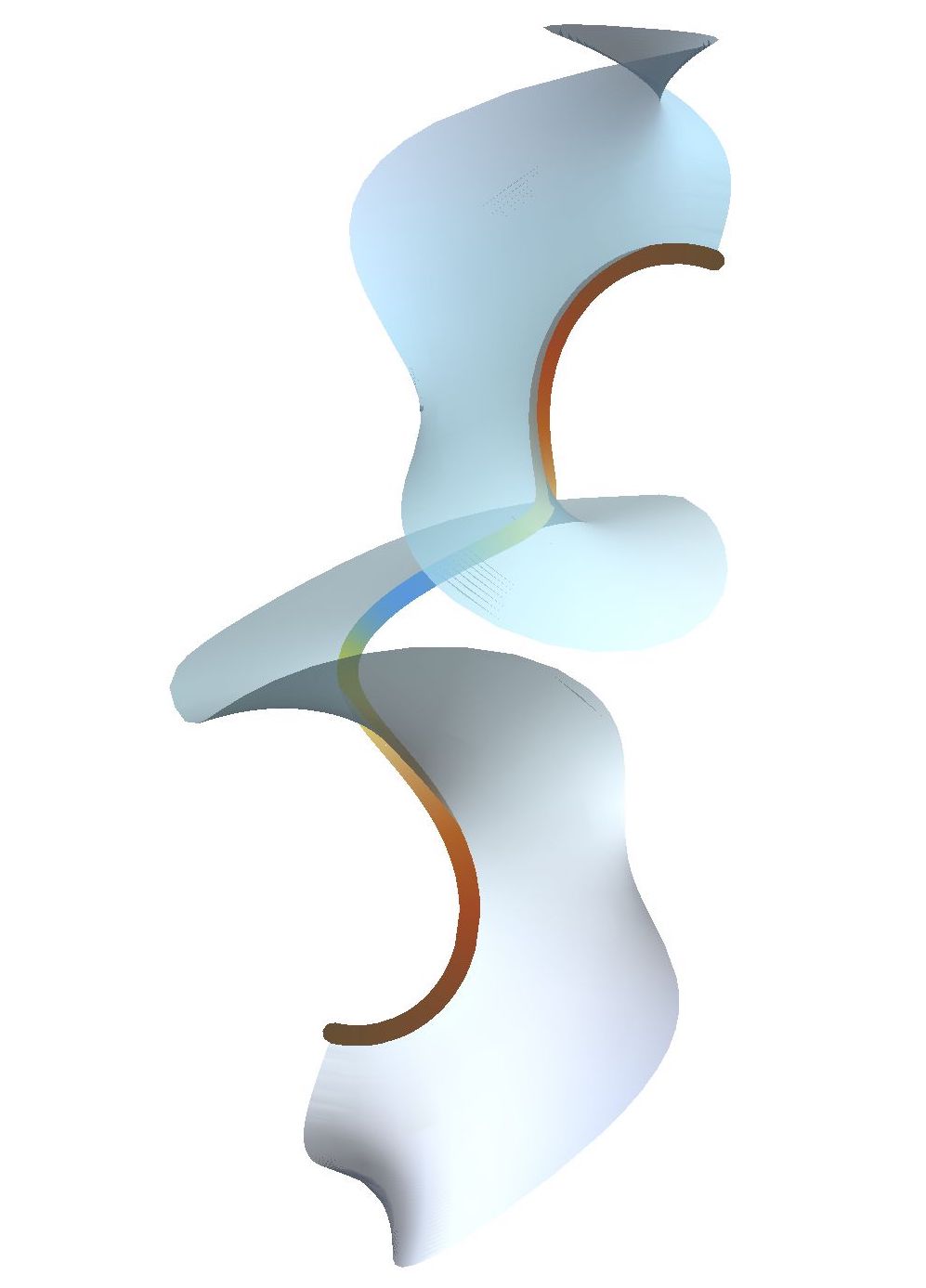}
\end{subfigure}
\quad
\begin{subfigure}[b]{0.3\linewidth}
\includegraphics[width=\linewidth]{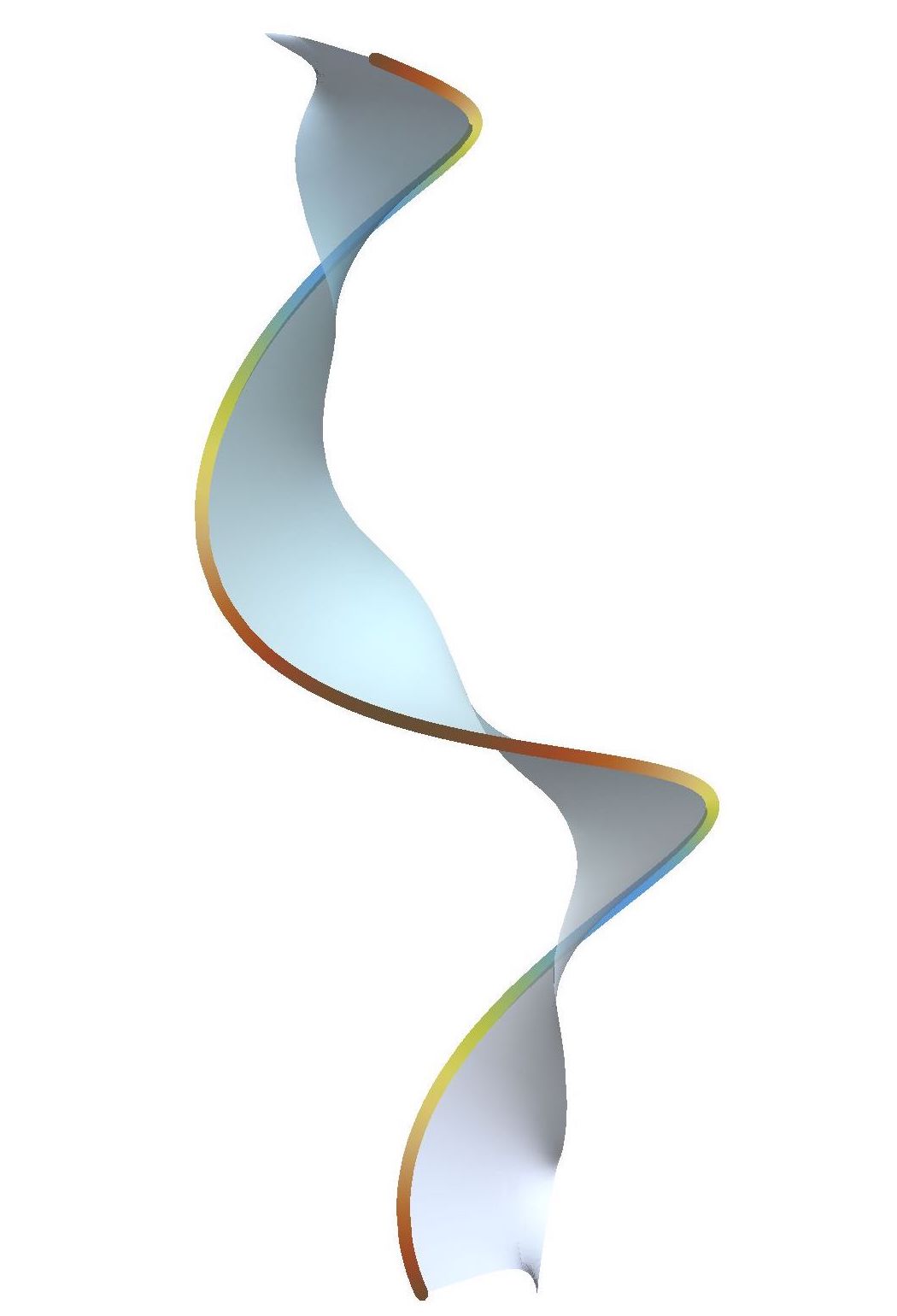}
\end{subfigure}
\caption{Three ``local" critical domains for $E\equiv E_{\sigma,\eta,\alpha,\beta}$ constructed in Mathematica using Bj\"orling's formula. From left to right: $E_{1,-5,1,1}$, $E_{1,1,-1,1}$ and $E_{1,6,6,0.11}$.}\label{Fig1}
\end{figure}

Next, let us consider the case where $\Sigma$ is a  surface of genus zero with an arbitrary number of boundary components, i.e. $\partial\Sigma\equiv\cup_{i=1}^m C_i$ with $m\geq 2$ (the case $m=1$ is covered in Theorem~\ref{disk}). In this setting, the following result holds as a particular case of Theorem \ref{thmplanar}. 

\begin{cor}\label{0} Let $X:\Sigma\rightarrow\mathbb{R}^3$ be an immersion of a genus zero surface critical for $E[X]$ defined in \eqref{energy} with $\eta\neq 0$ and suppose that $2\alpha\kappa_g+\eta\neq 0$ on at least one point of some boundary component $C\subset\partial\Sigma$. If $\tau_g=0$ holds somewhere on $C$, then the surface is a planar domain bounded by area-constrained elasticae.
\end{cor}

\begin{rem} For immersions of genus zero surfaces which satisfy the stronger hypothesis of $\tau_g\equiv 0$ holds along the entire boundary, the result of Corollary \ref{0} can also be proved using different techniques. Two essentially different proofs are sketched in Appendix B.
\end{rem}

On the other hand, assume that $\tau_g$ is not zero anywhere along $\partial\Sigma$. The following theorem shows that the topology of $\Sigma$ is prescribed by this condition.

\begin{thm}\label{ph} Let $\eta\neq 0$ and $X:\Sigma\rightarrow\mathbb{R}^3$ be an immersion of a genus zero surface critical for the energy $E[X]$, defined in \eqref{energy}. If $\tau_g>0$ (or $\tau_g<0$) everywhere along $\partial\Sigma$, then $\Sigma\cong A$ is a topological annulus.
\end{thm}
\textit{Proof.} Let $X:\Sigma\rightarrow\mathbb{R}^3$ be a critical immersion for $E[X]$. From the Euler-Lagrange equation \eqref{EL1}, $H\equiv 0$ holds on $\Sigma$, i.e. the surface is minimal. We then assume that $X:\Sigma\rightarrow\mathbb{R}^3$ is a conformal immersion of a bounded domain of the complex plane $\Sigma\subset\mathbb{C}$. From this, it follows that for an arbitrary complex coordinate $\omega$ the Hopf differential
$$\Phi\,d\omega^2=-\left(\mathbb{I}_{22}+i\mathbb{I}_{12}\right)d\omega^2,$$
is holomorphic in the bounded domain $\Sigma\subset\mathbb{C}$. Consequently, the imaginary part of $-\Phi$, i.e. $\mathbb{I}_{12}$, is harmonic.

Next, since $\eta\neq 0$, \eqref{EL2} implies that $\kappa_n\equiv 0$ holds along $\partial\Sigma$, so that $\Phi\lvert_{\partial\Sigma}=-i\tau_g e^\zeta = -\mathbb{I}_{12}\lvert_{\partial\Sigma}.$ The Minimum Principle for harmonic functions then implies that
$$\min_\Sigma \mathbb{I}_{12}=\min_{\partial\Sigma} \mathbb{I}_{12}=\min_{\partial\Sigma} \tau_g e^\zeta>0\,,$$
since $\tau_g> 0$ on $\partial\Sigma$.  (Equivalently, it follows from the Maximum Principle that $\max_\Sigma\mathbb{I}_{12}<0$ if  $\tau_g<0$ is assumed). Therefore, $\mathbb{I}_{12}$ is nonvanishing on $\Sigma$,  hence so is $\Phi$. Consequently, $\Sigma$ has no umbilic points.

Now, since the Hopf differential is nonvanishing on $\Sigma$, its horizontal foliation $\{ v\,|\, \Phi\,d\omega^2(v,v) = 0 \}$ is a global nonvanishing vector field on the surface.  Moreover, by considering the closed surface $2\Sigma$ generated by gluing two (appropriately oriented) copies of $\Sigma$ along their boundaries, this foliation also yields a global nonvanishing vector field on $2\Sigma$.  Therefore, it follows from the Poincar\'e-Hopf Index Theorem (see e.g. \cite{H}) that $\chi(2\Sigma) = 0$.   The relationship $\chi(2\Sigma) = 2\chi(\Sigma) - \chi(\partial \Sigma)$ then shows that the Euler-Poincar\'e characteristic of $\Sigma$ is zero, since $\partial\Sigma$ has odd dimension. Hence $\Sigma$ is a topological annulus, as claimed.\hfill$\square$
\\

Motivated by this result, it is interesting to study the annular case on its own. With this restriction, there is the following result. 

\begin{prop}\label{p} Let $\eta\neq 0$ and $X:\Sigma\cong A\rightarrow\mathbb{R}^3$ be an immersion of a topological annulus critical for $E[X]$, defined in \eqref{energy}. If $2\alpha\kappa_g+\eta\neq 0$ anywhere on $\partial\Sigma$, then precisely one of the following holds:
\begin{enumerate}
\item The domain $X(A)$ is planar and bounded by two area-constrained elasticae.
\item The geodesic torsion is everywhere positive (or negative) along the boundary, i.e. $\tau_g>0$ (resp., $\tau_g<0$) along $\partial A$.
\end{enumerate}
\end{prop}
\textit{Proof.} Let $X:A\rightarrow\mathbb{R}^3$ be critical for $E[X]$. We may assume the annulus $A$ is conformal to a domain in the complex plane $\mathbb{C}$ which is bounded by two circles $C_1$ and $C_2$, so that $\partial A\equiv C_1\cup C_2$ forms the (positively oriented) boundary. From the Euler-Lagrange equation \eqref{EL1}, the Hopf differential in the usual complex coordinate $z$, $\widetilde{\Phi}\,dz^2$, is holomorphic on $A$ and so the function $z^2\,\widetilde{\Phi}$ is analytic on $A$. Then, by Cauchy's Theorem 
\begin{equation*}
0=\oint_{\partial A} z^2\,\widetilde{\Phi}(z)\,dz=\oint_{C_1} z^2\,\widetilde{\Phi}(z)\,dz+\oint_{C_2} z^2\,\widetilde{\Phi}(z)\,dz=\int_{0}^{2\pi} \left[\tau_g e^\zeta\right]_{C_1}dt-\int_{0}^{2\pi}\left[\tau_g e^\zeta\right]_{C_2}dt\,,
\end{equation*}
where last equality holds after combining $\eta\neq 0$ with equation \eqref{EL2} and taking into account that along $\partial A$, $z^2\widetilde{\Phi}=\Phi=-i\tau_g e^\zeta$ holds, \eqref{Phiboundary}.

Now, suppose that there is at least one boundary point where $\tau_g=0$ holds. We may assume this point belongs to the boundary component $C_1$. In this case, the first conclusion follows immediately from Theorem \ref{thmplanar} (see also Corollary \ref{0}).  Otherwise, $\tau_g\neq 0$ everywhere on $\partial A$, and the above calculation shows that $\tau_g>0$ (resp. $\tau_g<0$) holds on $C_2$ whenever $\tau_g>0$ (resp. $\tau_g<0$) on $C_1$. This establishes the second conclusion. \hfill$\square$
\\

Annular domains appear naturally in the theory of minimal surfaces. Arguably, the most common way to characterize a minimal surface is via its Weierstrass representation \eqref{Weierstrass}, whose involved integrals may have periods on a non-simply connected domain. These immersions are generally multivalued, and for minimal surfaces $\Sigma$ it is common that a suitable quotient $\Sigma/\mathbb{Z}$ can be identified with an  annulus $A$. In this case, it follows that the first two fundamental forms determining the geometry of $\Sigma$ descend to $A$.

Lest the above results suggest that there are only planar equilibria of \eqref{energy}, we now construct a non-planar example. In particular, the goal is to find a multivalued immersion $X:A\rightarrow\mathbb{R}^3$ such that the image surface is critical for $E[X]$ in the quotient $\mathbb{R}^3/\mathbb{Z}$. Recall that the immersion $X:\mathbb{R}^2\rightarrow\mathbb{R}^3$ given by
\begin{equation}\label{helicoid}
X(r,\vartheta)=\left(r\cos\vartheta,r\sin\vartheta,a\vartheta+b\right)
\end{equation}
defines a minimal helicoid for given constants $a\neq 0$ and $b\in\mathbb{R}$. In this case, the curves corresponding to constant $r$ are helices, and along them the equations
$$\kappa_n\equiv 0\,,\quad\quad\quad \tau_g=\frac{-a}{a^2+r^2}\,,\quad\quad\quad \kappa_g=\frac{\pm r}{a^2+r^2}$$
hold. Now, since the surface is minimal, i.e. $H\equiv 0$, it is clear that the Euler-Lagrange equation \eqref{EL1} holds directly. Moreover, on helices the equation \eqref{EL2} also holds for any value $\eta\in\mathbb{R}$, since the normal curvature is identically zero. Finally, it is easy to check that for suitable constants $\sigma$, $\eta$, $\alpha$ and $\beta$, the boundary conditions \eqref{EL3}-\eqref{EL4} are also satisfied, since the above quantities are all constant. Consequently, the domains $\Omega$ in the helicoid defined by $r_0\leq r\leq r_1$ and $0\leq\vartheta\leq 2\pi$ for any constants $0<r_0<r_1$ correspond to minimal annuli in a quotient $\mathbb{R}^3/\mathbb{Z}$ which are critical for the energy $E[X]$, \eqref{energy}. 

Alternatively, these domains $\Omega$ can be understood as critical surfaces for $E[X]$, \eqref{energy}, having \emph{partially elastic boundary}, \cite{PP}. In this approach, the line segments $\vartheta=\vartheta_i$, $i=0,1$ are considered to be the fixed boundary components. Here, we should consider variations keeping these two segments fixed, i.e. $\delta X\equiv 0$ on $\vartheta_i$, $i=0,1$.

The helicoid has recently been used to model a stacked endoplasmic reticulum, which contributes to protein formation and transport in biological cells, \cite{Terasaki}. For this application, the multivalence of the immersion is an essential property, as it results in the stacking of membrane layers.  It is suspected that other minimizers of the Euler-Plateau energy with elastic modulus may have similar utility as models for biological phenomena.

\section*{Appendix A: Variation of the Total Geodesic Curvature}

The following calculation will show that the variation of the elastic modulus term $\int_\Sigma K\,d\Sigma$ on $\Sigma$ can be alternatively computed as the variation of the total geodesic curvature on $\partial\Sigma$. Note that the validity of this technique follows from the Gauss-Bonnet Theorem \eqref{GB}, which implies that these objects have identical variations.

To proceed, it is necessary to compute the pointwise variation of the geodesic curvature, $\kappa_g$. For this purpose, consider a general variation of $X:\Sigma\rightarrow\mathbb{R}^3$ whose restriction to the boundary is defined through $C(\epsilon) = C + \epsilon\,\delta C$.  In terms of the Darboux frame $\{T,\nu,n\}$, the variation $\delta C$ has the expression
\[ \delta C = \phi\,T+\psi\,\nu+\varphi\,n\,,\]
for some sufficiently smooth functions $\phi,\psi,\varphi$ on $\partial\Sigma$.  Now, using that $\delta T=\left[\left(\delta C\right)'\right]^\perp$ (here, $\left(\,\right)^\perp$ means orthogonal to $T$) and $\delta\nu=d\nu\left(\delta C^T\right)-\nabla\left(\nu\cdot \delta X\right)$, where $\left(\,\right)^T$ denotes the tangent component to the immersion $X$, a straightforward computation yields the variation of the Darboux frame with respect to $\delta C$,
\begin{eqnarray*}
\begin{cases}
T_\epsilon=\left(\psi'+\kappa_n\phi+\tau_g\varphi\right)\nu+\left(\varphi'+\kappa_g\phi-\tau_g\psi\right)n\,,\\
\nu_\epsilon=\left(\left[\kappa_n-2H\right]\varphi-\tau_g\phi-\partial_n\psi\right)n+\left(-\tau_g\varphi-\kappa_n\phi-\psi'\right)T\,,\\
n_\epsilon=\left(\tau_g\psi-\kappa_g\phi-\varphi'\right)T+\left(\partial_n\psi+\tau_g\phi+\left[2H-\kappa_n\right]\varphi\right)\nu\,,
\end{cases}
\end{eqnarray*}
where $\partial_n$ represents the derivative in the co-normal direction. Moreover, the geodesic curvature has the expression
$$\kappa_g=T'\cdot n=\frac{T_\rho}{\lVert C_\rho\rVert}\cdot n\,,$$
where $\rho$ denotes an arbitrary parameter and $F_\rho$ denotes the derivative of quantity $F$ with respect to $\rho$.  Therefore, using the variation of Darboux frame and differentiating above relation with respect to $\delta C$, a long but straightforward computation gives the desired pointwise variation,
$$\delta \kappa_g=\left(n\cdot \left[\delta C\right]'\right)'-\tau_g \left(\nu\cdot \delta C\right)'+\kappa_n\,\partial_n\left(\nu\cdot\delta X\right)+Kn\cdot\delta C-\kappa_g\, T\cdot\left(\delta C\right)'\,.$$
Finally, using the above together with integration by parts is sufficient for the expression
\begin{eqnarray*}
\delta\left(\oint_{\partial\Sigma}\kappa_g\,ds\right)&=&\oint_{\partial\Sigma}\delta\kappa_g\,ds+\kappa_g\,\delta\left(ds\right)=\oint_{\partial\Sigma}\left(-\tau_g \left(\nu\cdot \delta C\right)'+\kappa_n\,\partial_n\left(\nu\cdot\delta X\right)+Kn\cdot\delta C\right)ds\\&=&\oint_{\partial\Sigma}\left(\left[\tau_g'\,\nu+Kn\right]\cdot\delta C +\kappa_n\,\partial_n\left[\nu\cdot\delta X\right]\right)ds\,.
\end{eqnarray*} 

\section*{Appendix B: Alternative Proofs}
This Appendix presents some alternative arguments for selected results from the body.  First, recall that if a boundary component of a minimal surface is flat, then it is automatically planar (c.f Remark~\ref{plfl}). Moreover, it follows from this that the contact angle satisfies $\theta\equiv\pm\pi/2$. In this case, any of the following can be used to establish the conclusion of Proposition 3.8.

\begin{enumerate}
\item Assume that $C$ is the planar boundary component. Since the immersion is minimal, $C$ is composed entirely of umbilic points, which are also the zeros of the Hopf differential $\Phi\,d\omega^2$.  But, the Hopf differential is holomorphic on minimal surfaces (c.f. the proof of Theorem~\ref{disk}), so it must have isolated zeroes if it does not vanish identically.  Therefore, $\Phi\equiv 0$ and the surface is totally umbilical, hence planar.

\item Assume that the planar boundary component $C$ lies in a horizontal plane, and denote by $\mathcal{R}_t$ the one parameter family of rotations about any vertical axis. Then, the  function 
$$\psi:=\partial_t\left(\mathcal{R}_tX\right)_{t=0}\cdot \nu=E_3\times X\cdot \nu\,,$$
is the normal part of the derivative of a variation of $\Sigma$ through minimal surfaces.  As such, it follows that $\psi$ defines a Jacobi field along the surface. Moreover, the vector field $\nu$ is proportional to $E_3$ along $C$ (since $\kappa_n\equiv 0$), so $\psi\equiv 0$ also holds there. Using that $E_3$ is normal to $C$ along with $\tau_g\equiv 0$ and $H=\kappa_n=0$, it follows that 
$$\partial_n\psi=E_3\times n\cdot \nu+E_3\times X\cdot d\nu(n)=-T\cdot E_3-\tau_g E_3\times X\cdot T-\left(2H-\kappa_n\right)E_3\times X\cdot n=0\,.$$
Now, since $C$ is analytic it follows from the Cauchy-Kovalevskaya Theorem that the Cauchy problem
$$\Delta\psi+\lvert d\nu\rvert^2\psi=0,$$
with analytic initial conditions $\psi\equiv 0$ and $\partial_n\psi\equiv 0$ along $C$, has a unique (local) analytic solution $\psi\equiv 0$. Finally, using analyticity of the surface, it follows that $\psi\equiv 0$ must hold globally,  and so the surface $\Sigma$ is planar.

\item Assume that the planar boundary component $C$ lies in the horizontal plane $z=0$. Then, $\Sigma$ can be expressed locally as a graph with parameterization  $X(x,y)=\left(x,y,u(x,y)\right)$. With this, there is the Cauchy problem
$$\nabla\cdot\left(\frac{\nabla u}{\sqrt{1+\lvert \nabla u\rvert^2}}\right)=0\,,$$
with $u\equiv 0$ and $\nabla u\equiv 0$ along $C$. Note that the last condition comes from $\nu$ being vertical along $C$ (see \eqref{EL2}). Therefore, the Cauchy-Kovalevskaya Theorem shows that $u\equiv 0$ on a local domain containing $C$.  Moreover, it follows from the analyticity of $\Sigma$ that this solution must be global, and hence $u\equiv 0$ on the entirety of $\Sigma$. This proves that the surface is planar.
\end{enumerate}

In addition to this, there are the following alternative proofs of Corollary~4.2 under the (weaker) hypothesis that $\tau_g \equiv 0$ on $\partial\Sigma$:

\begin{enumerate}
\item Recall that the Hopf differential $\Phi\,d\omega^2$ is holomorphic on any surface satisfying \eqref{EL1}, e.g. on $\Sigma$. Moreover, from \eqref{Phiboundary} and the fact that $\tau_g\equiv 0$ along $\partial\Sigma$, it follows that $\Phi$ is defined on the bounded domain $\Sigma\subset\mathbb{C}$ and vanishes on the boundary $\partial\Sigma$. Therefore, the Maximum Principle implies that $\Phi\equiv 0$ holds on $\Sigma$, so (noticing that $H\equiv 0$ and the zeros of $\Phi$ are the umbilics) this shows that the surface is planar. Finally, using this information in \eqref{EL3} and \eqref{EL4} completes the argument.

\item Notice that the Gaussian curvature $K$ is nonpositive on minimal surfaces, i.e. $K\leq 0$. Therefore, using the Gauss-Bonnet Theorem \eqref{GB}, it follows that
$$0\geq \int_\Sigma K\,d\Sigma=\oint_{\partial\Sigma} \kappa_g\,ds+2\pi\chi(\Sigma)=2\pi\left(\left[m-1\right]-1\right)+2\pi\left(2-m\right)=0\,,$$
since the Euler-Poincar\'e characteristic of $\Sigma$ satisfies $\chi(\Sigma)=2-m$, where $m$ denotes the number of boundary components. Here it was also used that $\kappa_n\equiv 0$ holds along $\partial\Sigma$, so that $\kappa_g$ represents minus the signed curvature of the planar ($\tau_g\equiv 0$) boundary. Moreover, a version of the Gauss-Bonnet Theorem informally called as ``turning angles theorem" \cite{Milnor} was used to compute the total curvature of the boundary (observe that our choice of orientation coincides with that of Jordan). As a conclusion, $K\equiv 0$ identically on $\Sigma$ and the surface is planar.
\end{enumerate}

\bigskip

\begin{flushleft}
Anthony G{\footnotesize RUBER}\\
Department of Mathematics and Statistics, Texas Tech University-Costa Rica, San Jose, 10203, Costa Rica\\
E-mail: anthony.gruber@ttu.edu
\end{flushleft} 

\begin{flushleft}
\'Alvaro P{\footnotesize \'AMPANO}\\
Department of Mathematics and Statistics, Texas Tech University, Lubbock, TX, 79409, USA\\
E-mail: alvaro.pampano@ttu.edu
\end{flushleft} 

\begin{flushleft}
Magdalena T{\footnotesize ODA}\\
Department of Mathematics and Statistics, Texas Tech University, Lubbock, TX, 79409, USA\\
E-mail: magda.toda@ttu.edu
\end{flushleft} 

\end{document}